\documentclass{article}
\usepackage{amsfonts}

\usepackage{amsmath}
\usepackage{chicago}

\newtheorem{lemma}{Lemma}
\newtheorem{theorem}{Theorem}
\newtheorem{corollary}{Corollary}

\input{tcilatex}

\begin{document}

\title{Consistent Estimation of Pricing Kernels from Noisy Price Data}
\author{Vladislav Kargin\thanks{%
Cornerstone Research, 599 Lexington Avenue, New York, NY 10022, USA;
slava@bu.edu} }
\maketitle

\begin{abstract}
If pricing kernels are assumed non-negative then the inverse problem of
finding the pricing kernel is well-posed. The constrained least squares
method provides a consistent estimate of the pricing kernel. When the data
are limited, a new method is suggested: relaxed maximization of the relative
entropy. This estimator is also consistent.



MSC (2000) 62G08, 45Q05, 65J22
\end{abstract}

\bigskip \newpage

\section{Introduction}

Modern finance theory postulates that the price of a security is an integral
of its future payoff multiplied by a pricing kernel: 
\begin{equation}
S(x,\theta )=\int F(x^{\prime },\theta )p(x,x^{\prime })dx^{\prime }.
\end{equation}%
Here $S$ represents the security price, $x$ and $x^{\prime }$ current and
future values of stochastic factors relevant for pricing the security, $%
\theta $ a non-stochastic parameter, $F$ the future payoff, and $p$ the
pricing kernel. The pricing kernel is of great interest to finance theory
because it sheds light on investors' preferences over current and delayed
consumption. Practitioners are also interested in the pricing kernel because
it helps in pricing new securities, finding mispriced assets, and managing
risk\footnote{%
See, for example, applications in \shortciteN{jackwerth00}, %
\shortciteN{ait-sahalia_lo00}, and \shortciteN{rosenberg_engle02}.}.
Non-surprisingly, when it was discovered that the pricing kernel can be
recovered from option prices\footnote{%
By \shortciteN{breeden_litzenberger78} and \shortciteN{banz_miller78}, and
revived by \shortciteN{rubinstein94}.}, financial economists en masse went
agog inventing new and better methods for estimating the pricing kernel.%
\footnote{%
An incomplete list includes \shortciteN{jackwerth_rubinstein96}, %
\shortciteN{avellaneda97}, \shortciteN{soderlind_swensson97}, %
\shortciteN{melick_thomas97}, \shortciteN{ait-sahalia_lo98}, %
\shortciteN{avellaneda98}, \shortciteN{jackwerth00}, %
\shortciteN{ait-sahalia_duarte03}, and \shortciteN{bondarenko03}.} Many of
the methods, however, are heuristic and lack a rigorous proof of
consistency. This paper focuses on providing a simple proof of consistency
for the constrained least squares and a modified maximum entropy methods.

Mathematically, the pricing kernel estimation is an inverse problem. A
linear operator maps a set of functions (``pricing kernels'') into another
set of functions (``security prices''), and the problem is to invert this
operator. Often the problem is additionally complicated by the fact that
prices are observed only for a discrete set of securities and contaminated
with noise. This kind of inverse problems frequently appears in diverse
areas of applied mathematics and thoroughly studied \footnote{%
See reviews in \shortciteN{tikhonov_arsenin77}, \shortciteN{osullivan86},
and \shortciteN{engl00}.}.

The pricing kernel estimation is, however, special and what makes it special
is that the pricing kernel must be non-negative to prevent the existence of
systematic arbitrage opportunities.

This restriction on the operator's domain helps a lot. Without it, the
inverse problem is ill-posed, that is, the pricing operator does not have a
continuous inverse. Intuitively, small changes in security prices could lead
to large changes in the estimate of the pricing kernel. In addition, without
this restriction, the least squares method of estimation is inconsistent.
The pricing kernel selected by the least squares would fit the prices
exactly but would not converge to the true pricing kernel. In contrast,
non-negativity of the pricing kernel makes the corresponding inverse problem
well-posed and the least squares method consistent.

The key to the well-posedness is that for pricing purposes it is enough to
estimate the distribution function of the pricing kernel: cumulative pricing
kernel. These functions form a Banach space with respect to the uniform
convergence topology, and we can apply one of the Banach theorems: A
continuous one-to-one operator on a Banach space has a continuous inverse.
Consequently, the inverse problem is well-posed.

What can be said about consistency? By well-posedness, estimating the kernel
can be reduced to estimating price function from noisy observations: the map
from price functions to pricing kernels is continuous and cannot inflate the
error of estimation. Luckily, the problem of estimating the price function
is the classic problem of non-parametric estimation of a regression
function, and for this problem the conditions of the least squares
consistency are well known. It turns out that they are satisfied provided
the pricing kernel is non-negative. Intuitively, additional information
about the structure of pricing kernels prevents overfitting of the
regression function and forces consistent convergence of the estimates.
Together with well-posedness, this implies that the constrained least
squares estimates the cumulative pricing kernel consistently.

While asymptotically consistent, the constrained least squares may, however,
perform unsatisfactorily in small samples. It is because this method ignores
prior information about the pricing kernel. One way to abate the problem is
to include in the objective function a term that measures distance from the
prior information pricing kernel. This idea leads to a method that combines
advantages of both the least squares and the maximal entropy methods. The
method maximizes the weighted sum of relative entropy and the mean squared
pricing error. With suitably chosen parameters, this method is also
consistent.

Let me briefly describe the related literature. The maximum entropy method
for estimating pricing kernel was developed by \shortciteN{buchen_kelly96}
and \shortciteN{stutzer96}, following a suggestion in %
\shortciteN{rubinstein94}, and elaborated by \shortciteN{avellaneda97}, %
\shortciteN{avellaneda98}, and \shortciteN{fritelli00}. These papers
typically assume that the securities are priced correctly but only a scarce
discrete set of prices is known. For the alternative case of large amount of
noisy data, methods of pricing kernel estimation based on smoothing or other
ideas were developed by \shortciteN{jackwerth_rubinstein96}, %
\shortciteN{ait-sahalia_lo98}, \shortciteN{jackwerth00}, and %
\shortciteN{bondarenko03} among others. Implicitly, these papers address the
problem of ill-posedness of kernel estimation by the classic method of
regularization. This paper is different because it shows that on the
restricted domain of non-negative kernels the problem is well-posed and so
does not need additional regularization.

\shortciteN{ait-sahalia_duarte03} estimate the pricing kernels by smoothing
the \emph{constrained} least squares estimator, and refer to the statistical
literature for the proof of consistency. In this paper, we provide an
explicit proof of the constrained least squares consistency and consider
another modification of the method based on \ the idea of entropy distance
minimization.

The rest of the paper is organized as follows. Section 2 reminds the basics
of the theory of pricing kernels. Section 3 shows that the problem of
finding the non-decreasing cumulative pricing kernel is well-posed. Section
4 demonstrates consistency of the least squares. Section 5 explains how the
idea of the maximum entropy can be used to improve the least squares method,
and proves the consistency of the modification. Section 6 concludes.

\section{What is the pricing kernel?}

The pricing kernel, $p$, is a function of stochastic factors that allows
pricing securities by using their future payoff functions:%
\begin{equation}
S(x,\theta )=\int F(x^{\prime },\theta )p(x,x^{\prime })dx^{\prime }.
\end{equation}%
Often, the choice of units in which the stochastic factors are measured is
arbitrary, so we can normalize the initial level of factors: $x=1.$ Slightly
abusing notation, we will denote $p(1,x)$ as $p(x).$ Let us define \emph{%
cumulative pricing kernel} as follows$:$%
\begin{equation}
P(x)=\int_{-\infty }^{x}p(t)dt.
\end{equation}
With these notations, the pricing formula can be rewritten in a more
convenient form:%
\begin{equation}
S(\theta )=\int F(x,\theta )dP(x).
\end{equation}

Non-negativity of the pricing kernel, implied by the absence of arbitrage
opportunities (\shortciteN{harrison_kreps79})$,$ translates into
monotonicity of the cumulative pricing kernel: $P(x)$ is non-decreasing. In
addition, the price of the security that have a unit payoff is finite, so $%
P(x)$ is bounded.

We will be interested in pricing kernels that depend only on one factor. For
example, when the class of securities consists of options written on another
security, this factor is the price of the underlying security. We can
further simplify the problem by noting that payoff of most derivative
securities that occur in practice can be represented as a linear combination
of underlying security and security that has a non-zero payoff only if the
price of underlying is less than a certain bound, $B$. Therefore, we can
concentrate on pricing the derivatives with finite support, and then by
integration by parts we have the following pricing formula:%
\begin{equation}
S(\theta )=-\int_{0}^{B}P(x)dF(x,\theta ),  \label{pricing_formula}
\end{equation}%
where $B$ is such that $F(x,\theta )=0$ for $x>B.$

The next lemma shows that estimating cumulative pricing kernel is sufficient
for pricing purposes. Let $P_{n}(x)$ be an estimate of $P(x)$.\ Let $S_{n}$
be the corresponding price of the derivative from (\ref{pricing_formula}).

\begin{lemma}
\label{continuity}Suppose $F(x)$ has bounded variation and $P_{n}$ converges
to $P$ in uniform metric as $n$ goes to $\infty $. Then $S_{n}$ converges to 
$S.$
\end{lemma}

\textbf{Proof:} 
\begin{eqnarray}
\left| S_{n}-S\right| &=&\left| \int_{0}^{B}\left[ P_{n}(x)-P(x)\right]
dF(x)\right| \leq \int_{0}^{B}\left| P_{n}(x)-P(x)\right| dF(x) \\
&\leq &C\left\| P_{n}(x)-P(x)\right\| _{\infty },
\end{eqnarray}%
where $C$ is the total variation of $F(x).$ QED.

Consider now how we can estimate the pricing kernel. Typically, it is done
by using the prices of puts. A European put with strike $K$ is a security
that will pay:%
\begin{equation}
F(x,K)=\max \left( K-x,0\right)
\end{equation}%
at the expiration date if the price of the underlying security is $x$ on
that date. Then, the price of the put with strike $K$ is 
\begin{eqnarray}
S(K) &=&\int_{0}^{K}\left( K-x\right) dP(x) \\
&=&\int_{0}^{K}P(x)dx.
\end{eqnarray}%
The operator of interest is then 
\begin{equation}
A:P(x)\rightarrow S(K)=\int_{0}^{K}P(x)dx.  \label{pricing_operator}
\end{equation}%
In a more general setting, we are interested in the inverse problem defined
by operator%
\begin{equation}
A_{F}:P(x)\rightarrow S(\theta )=-\int_{0}^{B}P(x)dF(x,\theta ).
\label{pricing_operator2}
\end{equation}%
As Lemma \ref{continuity} shows, operator $A_{F}$ is continuous in the
uniform metric ($L^{\infty }$). We are interested in knowing whether its
inverse is continuous, that is, if small deviations in prices can lead to
large deviations in pricing kernel. We also need to know if the pricing
kernel can be consistently estimated from noisy and discrete data. These
problems are handled in the next two sections.

\section{Pricing problem is well-posed.}

A problem $Ax=y$ is called well-posed if the operator $A$ has a continuous
inverse. It is implicit in this definition that the operator is given with
its domain, and that topologies in both the range and the domain are
specified: the same operator may be ill-posed on one domain and well-posed
on another one. The concept of well-posedness originated in mathematical
physics by Hadamard as a tool to select the linear problems that could arise
from a physical problem. Later, however, it was discovered that many
important problems are ill-posed and the methods of their solution were
derived (\shortciteN{tikhonov_arsenin77}, \shortciteN{osullivan86}, %
\shortciteN{engl00}).

If no restrictions on pricing kernels were imposed, then the operator in (%
\ref{pricing_operator}) would correspond to an ill-posed problem. Indeed, it
is easy to see ill-posedness from the following example:%
\begin{equation}
\alpha \cos (\beta x)\rightarrow \int_{0}^{K}\alpha \cos (\beta x)dx=\frac{%
\alpha }{\beta }\sin (\beta K).
\end{equation}%
Consider the uniform convergence metric on both the domain and the range of
the operator. If we set $\beta =\alpha /\varepsilon ,$ then the norm of the
function on the left-hand side is constant: $\left\| \alpha \cos (\beta
x)\right\| _{\infty }=\alpha ,$ but its image can be made arbitrarily close
to zero $\left\| F(\alpha \cos (\beta x))\right\| _{\infty }=\left\|
\varepsilon \sin (\beta K)\right\| _{\infty }\leq \epsilon :$ therefore the
inversion operator acts discontinuously.

However, for the restricted domain of cumulative pricing kernels we have the
following theorem:

\begin{theorem}
\label{well_posedness} If $A_{F}$ is injective then it defines a well-posed
problem on the space of all non-decreasing continuous functions with uniform
convergence topology.
\end{theorem}

\textbf{Proof: } In uniform convergence topology the space of non-decreasing
continuous functions is complete. If $A_{F}$ is injective, then it defines a
continuous one-to-one correspondence between this space and its image. The
conclusion of the theorem follows because of one of the Banach theorems (see
for example Theorem 11 in Chapter 15 of \shortciteN{lax02}): A linear
operator that establishes a continuous one-to-one correspondence between two
complete normed linear spaces has a continuous inverse. QED.

\begin{corollary}
\label{well_posedness2}Operator $A$ from (\ref{pricing_operator}) defines a
well-posed problem on the space of all non-decreasing continuous functions
with uniform convergence topology.
\end{corollary}

\textbf{Proof:} Since $A$ is injective, Theorem \ref{well_posedness} can be
applied. QED.

In practice securities prices are known up to an error. This error includes
bid-ask spread, non-stationarity in the pricing kernel, market
inefficiencies and so on. We will use Theorem \ref{well_posedness} and
Corollary \ref{well_posedness2} as tools to prove that as the amount of data
grows the constrained least squares estimates the pricing kernel
consistently.

\section{Estimation by least squares is consistent.}

Let $\left\{ \Omega ,\Sigma ,\Pr \right\} $ be a probability space and $%
\varepsilon _{i}$ be a sequence of independent identically distributed
random variables with zero expectation and finite variance. Let also $x_{i}$
be a sequence of points located between $0$ and $B,$ which has a positive
density on $[0,B]$. We will say that the constrained least squares estimates
function $f$ from set $\mathcal{F}$ consistently in norm $\left\| \cdot
\right\| $ relative to operator $A$ if for any $\delta ,$ with probability $%
1 $ there exists such $N_{0}$ that for $N\geq N_{0},$ there exists 
\begin{equation}
f_{N}=\arg \max_{\widehat{f}\in \mathcal{F}}\sum_{i=1}^{N}\left( Af(x_{i})-A%
\widehat{f}(x_{i})\right) ^{2},
\end{equation}%
and $\left\| f-f_{N}\right\| <\delta .$ In other words, with probability 1$,$
the sequence of estimates $f_{N}$ converges to the true function.

Here we are interested in the set, $\mathcal{P}$, of non-decreasing,
continuous, bounded functions on interval $[0,B].$ We use the uniform
convergence topology and operator $A$ from (\ref{pricing_operator})$.$

\begin{theorem}
\label{consistency}The constrained least squares estimates any function in $%
\mathcal{P}$ consistently in $L^{\infty }$ relative to operator $A$.
\end{theorem}

\textbf{Proof: \ }Let $\mathcal{R}$ be the image of $\mathcal{P}$ under
operator $A$. Then $\mathcal{R}$ is the set of convex, increasing,
continuous, bounded functions. Because of Theorem \ref{well_posedness},
operator $A$ has a continuous inverse from $\mathcal{R}$ to $\mathcal{P}.$
Consequently, consistency of estimating functions from $\mathcal{P}$
relative to operator $A$ is equivalent to consistency of estimating
unmodified functions from $\mathcal{R}.$ For $\mathcal{R}$, we can apply
classic results for the non-parametric estimation of convex function. In
particular, according to the main Theorem in \shortciteN{hanson_pledger76},
convex functions can be estimated consistently in $L^{\infty }$ norm by the
constrained least squares method. QED.

For a more general operator $A_{F}$ from (\ref{pricing_operator2}) we have a
similar theorem, which, however, needs a more advanced technique and comes
to a weaker conclusion.

Let us call payoff function $F(x,\theta )$ \emph{uniformly Lipschitz} in $%
\theta $ if 
\begin{equation}
\left| F(x,\theta _{1})-F(x,\theta _{2})\right| \leq C\left| \theta
_{1}-\theta _{2}\right| ,
\end{equation}%
where $C$ does not depend on $x.$ Also let us call $F(x,\theta )$ \emph{%
uniformly bounded in variation}, if its total variation over $x\in \lbrack
0,B]$ is bounded by a constant that does not depend on $\theta .$

\begin{theorem}
\label{consistency2}If $F(x,\theta )$ is uniformly Lipschitz in $\theta $
and uniformly bounded in variation, and $A_{F}$ is injective, then the
constrained least squares estimates any function in $\mathcal{P}$
consistently in $L^{2}$ relative to operator $A_{F}$.
\end{theorem}

In the proof we will again aim to prove that any function in $\mathcal{R}%
=A_{F}(\mathcal{P})$ can be estimated consistently by the constrained least
squares. We are going to do it by referring to a theorem in %
\shortciteN{geer87}. First, let us introduce several additional concepts.
Let $X$ be a set of functions on $\mathbb{R}^{k}$ and let $M_{n}(\delta ,%
\mathcal{X})$ be the minimal number of elements in a $\delta -$covering of
set $\mathcal{X}$, if the distance is measured by the norm 
\begin{equation}
\left\| f\right\| _{n}=\frac{1}{n}\sum_{i=1}^{n}\left[ f(x_{i})\right] ^{2}.
\label{nnorm}
\end{equation}%
Then $\delta -$\emph{entropy} of a set is defined as 
\begin{equation}
N_{n}(\delta ,\mathcal{X})=\frac{1}{n}\log M_{n}(\delta ,\mathcal{X}).
\end{equation}%
Note that $\delta -$\emph{entropy }depends on the choice of points $x_{i}.$
We assume that they are distributed randomly according to a measure, $\mu ,$
that has a positive continuous density on $[0,B].$ Then let us call a set of
functions \emph{entropically thin} if for any $\delta $ 
\begin{equation}
N_{n}(\delta ,\mathcal{X})\rightarrow _{\mu }0\text{ as }n\rightarrow \infty
,
\end{equation}%
where convergence is in probability. Intuitively, a set of functions is
entropically thin if all its functions can be well approximated by functions
from a relatively ``small''\ subset.\footnote{%
The concept of entropy in relation to totally bounded sets of functions was
introduced by \shortciteN{kolmogorov_tikhomirov59}. It \ was applied to the
problem of consistency in non-parametric estimation by %
\shortciteN{vapnik_cervonenkis81}. For a textbook presentation, see %
\shortciteN{pollard84}.}

Next, a class of functions, $\mathcal{X},$ is called \emph{uniformly square
integrable} if 
\begin{equation}
\lim_{C\rightarrow \infty }\sup_{f\in \mathcal{X}}\int_{\left| f\right|
>C}f^{2}dx=0.
\end{equation}%
A somewhat weaker version of van de Geer's result is sufficient for our
purposes. It says that if a set $\mathcal{X}$ is uniformly square integrable
and entropically thin, then the constrained least squares method is $L^{2}-$%
consistent.

\textbf{Proof of Theorem \ref{consistency2}: }Any $S\in \mathcal{R}$ is
representable as 
\begin{equation}
S(\theta )=-\int_{0}^{B}P(x)dF(x,\theta ).
\end{equation}%
Since set $\mathcal{P}$ is uniformly bounded and $F(x,\theta )$ is uniformly
bounded in variation, set $\mathcal{R}$ is also uniformly bounded.
Consequently, it is uniformly square integrable.

Similarly, since $F(x,\theta )$ is uniformly Lipschitz in $\theta ,$ and $%
\mathcal{P}$ is uniformly bounded$,$ set $\mathcal{R}$ is uniformly
Lipschitz:%
\begin{eqnarray}
\left| S(\theta _{1})-S(\theta _{2})\right| &=&\left| \int_{0}^{B}\left[
F(x,\theta _{1})-F(x,\theta _{2})\right] dP(x)\right| \\
&\leq &\int_{0}^{B}C_{1}\left| \theta _{1}-\theta _{2}\right| dP(x) \\
&\leq &BC_{1}C_{2}\left| \theta _{1}-\theta _{2}\right| .
\end{eqnarray}
Consequently, by Lemma 3.3.1 in \shortciteN{geer87} $\mathcal{R}$ is
entropically thin. Therefore, van de Geer's theorem can be applied and the
constrained least squares estimator is $L^{2}-$consistent.

QED.

The conditions of Theorem \ref{consistency2} are not very restrictive. For
example, the set of payoff functions for puts is uniformly Lipschitz:%
\begin{equation}
\left| \max \left\{ K_{1}-x,0\right\} -\max \left\{ K_{2}-x,0\right\}
\right| \leq \left| K_{1}-K_{2}\right| .
\end{equation}%
It is also clearly uniformly bounded in variation over $x\in \lbrack 0,B],$
provided that we consider only a bounded set of the strikes: $\max \left\{
K-x,0\right\} \leq \overline{K}\equiv \max K.$ Finally, the pricing
operator, $A,$ is injective if $\overline{K}\geq B.$ Therefore, Theorem \ref%
{consistency2} is applicable.

While asymptotically consistent, the constrained least squares may perform
poorly in small samples. It fails to take into account such possible prior
beliefs as that the pricing kernel is smooth, or unimodal, or that it is
approximately proportional to an infinitely divisible probability
distribution, etc. In the next section we consider a modification of the
method of constrained least squares that allows to take into account the
prior information.

\section{Relaxed Maximum Relative Entropy Method}

In this section we will for simplicity restrict the discussion to the case
when the pricing kernel is estimated from the put prices. Relaxed maximum
entropy method penalizes both the degree to which the model fails in
explaining the price data and the model's deviation from a prior model:%
\begin{equation}
\widehat{P}(x)=\arg \min_{P(x)}\left\{ \frac{1}{N}\sum_{i=1}^{N}\left(
S(K_{i})-S_{i}\right) ^{2}+\lambda _{N}\int_{0}^{B}\ln \frac{dP(x)}{dP_{0}(x)%
}dP(x)\right\} ,
\end{equation}%
where $S_{i}$ is the observed price of the put with strike $K_{i}$, 
\begin{equation}
S(K)=AP(x)\equiv \int_{0}^{K}P(x)dx,
\end{equation}%
and $P_{0}(x)$ is a prior cumulative pricing kernel.

Recall that the regular maximum entropy method is described by the following
minimization problem:%
\begin{equation}
\widehat{P}_{ME}(x)=\arg \min_{P(x)}\left\{ \int_{0}^{B}\ln \frac{dP(x)}{%
dP_{0}(x)}dP(x)\text{ s.t. }S(K_{i})=S_{i}\text{ for each }i\right\} .
\end{equation}%
If prices are contaminated with noise, then the regular maximum entropy may
run into difficulties with the existence of the solution and is unlikely to
be consistent. In the relaxed maximum entropy method, constraints are not
rigid, they are substituted with a penalizing term in the objective
function. Consequently, the solution is guaranteed to exist. What about
consistency?

\begin{theorem}
\label{RME_consistency}There is such a sequence of positive constants $%
\lambda _{N},$ that the relaxed maximum entropy method estimates the
cumulative pricing kernel, $P(x),$ consistently in $L^{2}$ norm.
\end{theorem}

\textbf{Proof: }By a lemma below, there is such a sequence $\lambda
_{N}\rightarrow 0$ that as $N\rightarrow \infty ,$ the solution of the
problem\textbf{\ }%
\begin{equation}
\min_{S(K)}\left\{ \frac{1}{N}\sum_{i=1}^{N}\left( S(K_{i})-S_{i}\right)
^{2}+\lambda _{N}\int_{0}^{B}\ln \frac{dP(x)}{dP_{0}(x)}dP(x)\right\}
\end{equation}%
with probability 1 approaches in $L^{2}$ norm the solution of the
constrained least squares problem: 
\begin{equation}
\min_{S(K)}\frac{1}{N}\sum_{i=1}^{N}\left( S(K_{i})-S_{i}\right) ^{2}\text{
s.t. }\partial _{K}^{2}S\geq 0.
\end{equation}%
by Theorem \ref{consistency}, as $N\rightarrow \infty ,$ the solution of the
constrained least squares problem with probability 1 approaches the true
pricing function $S(K)$.

By the standard diagonal process argument there is such a sequence of $%
\lambda _{N},$ that the solution of 
\begin{equation}
\min_{S(K)}\left\{ \frac{1}{N}\sum_{i=1}^{N}\left( S(K_{i})-S_{i}\right)
^{2}+\lambda _{N}\int_{0}^{B}\ln \frac{dP(x)}{dP_{0}(x)}dP(x)\right\}
\end{equation}%
approaches in $L^{2}$ the true function $S(K)$ as $N\rightarrow \infty .$
Because the differentiation is a continuous operator on the set of convex
non-decreasing functions, $P(x)$ is also estimated consistently in $L^{2}.$
QED.

In the proof of Theorem \ref{RME_consistency}, we have used the following
Lemma. Consider the problem: 
\begin{equation*}
\min_{R\in \mathcal{R}}\left\{ F_{N}(R)+\lambda _{N}G(R)\right\} ,
\end{equation*}%
where $F_{N}$ and $G$ are continuous functionals of $R(x).$ Let $R_{N}$ be
solution of the problem,$\ $and $\widehat{R}_{N}$ be the solution for $%
\lambda =0.$ Let the sequence of functionals $F_{N}$ be called \emph{proper }%
on $\mathcal{R}$ if for any $\varepsilon $ we can find such $\delta $ that
for all sufficiently large $N,$ and $R\in \mathcal{R}$, $\ $condition $%
F_{N}(R)-F_{N}(\widehat{R}_{N})<\delta $ implies that $\left\Vert R-\widehat{%
R}_{N}\right\Vert _{L^{2}}\leq \varepsilon .$

\begin{lemma}
If $\{F_{N}\}$ is proper on $\mathcal{R}$, then there exists such a sequence 
$\lambda _{N}$ that $R_{N}$ $-\widehat{R}_{N}$ converges to zero.
\end{lemma}

\textbf{Proof: }Take an $\varepsilon $ and select $\delta $ and $N_{0}$ as
in the definition of properness; then for any $R\in \mathcal{R}$ and any $%
N\geq N_{0}$ from $\left\Vert R-\widehat{R}_{N}\right\Vert
_{L^{2}}>\varepsilon $ it follows that $F_{N}(R)-F_{N}(\widehat{R}_{N})\geq
\delta .$ On the other hand, from continuity of $F_{N}$ it follows that we
can find such $\varepsilon _{1}$ that $\left\Vert R-\widehat{R}%
_{N}\right\Vert _{L^{2}}<\varepsilon _{1}$ implies $F_{N}(R)-F_{N}(\widehat{R%
}_{N})\leq \delta /2.$ Also, since $G$ is continuous, we can find such an $%
R^{\prime }$ inside the $\varepsilon _{1}-$neighborhood of $\widehat{R}_{N}$
that $|G(R^{\prime })-G(\widehat{R}_{N})|<c.$ Consequently we can find such $%
\lambda $ that $\lambda G(R^{\prime })<\delta /2.$ Then it is clear that the
maximizer of $F_{N}(R)+\lambda G(R)$ cannot be outside of the $\varepsilon -$%
neighborhood of $\widehat{R}_{N}:$ $R^{\prime }$ would improve on it.

Thus, for any $\varepsilon ,$ there is $N_{0}$ and $\lambda $ such that for $%
N\geq N_{0}$ the solution of $\min \{F_{N}+\lambda G\}$ is in $\varepsilon -$%
neighborhood of $\widehat{R}_{N}.$ QED.

This Lemma can be used the proof of Theorem \ref{RME_consistency} because
the functional 
\begin{equation}
F_{N}(R)\equiv \frac{1}{N}\sum_{i=1}^{N}\left( R(K_{i})-S_{i}\right)
\end{equation}%
is proper. Indeed, note that this functional has a nice special property:

\begin{lemma}
If $F_{N}(R)-F_{N}(\widehat{R}_{N})<\varepsilon $ then $F_{N}(R-\widehat{R}%
_{N})<\varepsilon .$
\end{lemma}

\textbf{Proof:} Let $R=\widehat{R}_{N}+\delta R.$ Since $F_{N}$ is a
quadratic form, we can define a corresponding bilinear product:%
\begin{equation}
(f,g)=\frac{1}{2}\left\{ F_{N}(f+g)-F_{N}(f)-F_{N}(g)\right\} ,
\end{equation}%
We claim that 
\begin{equation}
(\widehat{R}_{N},\delta R)\geq 0.  \label{tightness_proof1}
\end{equation}%
Indeed, since the set of convex non-decreasing functions, $\mathcal{R},$ is
convex, $R_{\alpha }\equiv \widehat{R}_{N}+\alpha \delta R\in \mathcal{R}$
for any $\alpha \in \lbrack 0,1]$. Consequently, if (\ref{tightness_proof1})
were violated we could find such $\alpha $ that $(R_{\alpha },R_{\alpha })<(%
\widehat{R}_{N},\widehat{R}_{N}),$ which would contradict optimality of $%
\widehat{R}_{N}.$ Using (\ref{tightness_proof1}), \ we can write:%
\begin{equation}
F_{N}(\delta R)+F_{N}(\widehat{R}_{N})\leq F_{N}(R)<F_{N}(\widehat{R}%
_{N})+\varepsilon ,
\end{equation}%
and $F_{N}(\delta R)<\varepsilon .$ QED.

So, to obtain properness of $\{F_{N}\}$ it remains to prove that from $%
F_{N}(R-\widehat{R}_{N})<\varepsilon $ for all large $N$ we can conclude
that $\left\Vert R-\widehat{R}_{N}\right\Vert _{L_{2}}\leq \varepsilon .$ By
assumption, points $\{x_{i}\}$ are distributed with density $\rho (x)\geq
k>0 $ on $[0,B].$ Then we can use the following lemma:

\begin{lemma}
If $f(x)$ is non-negative and has finite variation on $[0,B]$ then from%
\begin{equation}
\frac{1}{N}\sum_{i=1}^{N}f(x_{i})\leq \varepsilon \text{ for any }N\geq
N_{0},
\end{equation}%
it follows that 
\begin{equation}
\int_{0}^{B}f(x)dx\leq \frac{\varepsilon }{k}.
\end{equation}
\end{lemma}

\textbf{Proof:} The sum converges to $\int_{0}^{B}f(x)\rho (x)dx.$ Therefore,%
\begin{equation}
\int_{0}^{B}f(x)dx=\int_{0}^{B}\frac{f(x)}{\rho (x)}\rho (x)dx\leq \frac{1}{k%
}\int_{0}^{B}f(x)\rho (x)dx\leq \frac{\varepsilon }{k}.
\end{equation}%
QED.

Properness of functional $F_{N}(x)$ follows from this lemma applied to $%
f(x)=\left( R(x)-\widehat{R}_{N}(x)\right) ^{2}.$

\section{Conclusion}

It is proved that the mapping from the set of non-decreasing cumulative
pricing kernels to security prices corresponds to a well-posed inverse
problem, and that the constrained least squares method provides a consistent
estimator of the cumulative pricing kernel.

It is also suggested that in small samples the performance of the
constrained least squares can be improved by a modification that takes into
account that the pricing kernel should be close to a certain prior kernel.
It is proved that this method is consistent.

\bibliographystyle{CHICAGO}
\bibliography{comtest}

\end{document}